\documentclass[12pt]{article}
\usepackage{amsfonts}
\usepackage{latexsym}
\usepackage{amsmath}

\newcommand{\zz}{\mathbb{Z}}
\newcommand{\pp}{\mathbb{P}}
\newcommand{\rr}{\mathbb{R}}
\newcommand{\cc}{\mathbb{C}}
\newcommand{\qq}{\mathbb{Q}}
\newcommand{\ff}{\mathbb{F}}
\newcommand{\lan}{\langle}
\newcommand{\ran}{\rangle}
\newcommand{\B}{\mathcal{B}}
\newcommand{\Bi}{\mathcal{B}^{-1}}

\newcommand{\Del}{\Delta^{\circ}}

\newcommand{\Ns}{N_{\sigma}}
\newcommand{\Gs}{G_{\sigma}}
\newcommand{\Gt}{G_{\tau}}
\newcommand{\mvi}{\lan m,v_i \ran}
\newcommand{\mvk}{\lan v_k,m \ran}
\newcommand{\msvj}{\lan m^*,v_j \ran}

\newcommand{\Us}{U_{\sigma}}
\newcommand{\Usp}{U_{\sigma'}}
\newcommand{\Ut}{U_{\tau}}
\newcommand{\Ms}{M_{\sigma}}

\newcommand{\csig}{\check{\sigma}}
\newcommand{\ssig}{S_{\sigma}}
\newcommand{\stau}{S_{\tau}}
\newcommand{\csigp}{\check{\sigma'}}

\newtheorem{theorem}{Theorem}

\newtheorem{corollary}{Corollary}

\newtheorem{lemma}{Lemma}

\begin{document}
\title{Orbifold Hodge Numbers of Calabi-Yau Hypersurfaces}
\author{Mainak Poddar\\
     University of Wisconsin-Madison\\
     Email: poddar@math.wisc.edu}
\date{}     
\maketitle
\begin{abstract}
We identify the twisted sectors of a compact simplicial toric variety.
We do the same for  a generic nondegenerate Calabi-Yau hypersurface of
an $n$-dimensional simplicial Fano toric variety and then explicitly 
compute $h^{1,1}_{orb}$ and $h^{n-2,1}_{orb}$ for the hypersurface. We 
give applications to the orbifold string theory conjecture and orbifold 
mirror symmetry.
\end{abstract}

\tableofcontents

\section{Introduction}

\noindent Recently, there is an emerging  subject of mathematics: 
 \emph{Stringy geometry and topology of orbifolds}. The core of this
 includes the developing theories of orbifold cohomology, orbifold 
 quantum cohomology, orbifold K-theory, motivic integration, Mckay
 correspondence etc.\\ 
   
\noindent  Physicists believe that orbifold string theory is equivalent
 to ordinary string theory of certain desingularisations. This belief  
 motivated a body of conjectures, collectively referred to as the
 \emph{Orbifold string theory conjecture}. The conjecture we are 
 interested in is the \emph{K-Orbifold string theory conjecture}.  
 It states that there is a natural isomorphism between the 
 \emph{Orbifold K-theory} of a Gorenstein orbifold and the ordinary
 K-theory of its crepant resolution.\\

\noindent To construct a natural isomorphism as the conjecture demands,
 is a very hard problem.
 But many weaker versions of the conjecture that compare Euler numbers, 
 Hodge numbers, etc. has been studied extensively in the literature in
 the case of orbifolds that are global-quotients.  Batyrev[2], and Batyrev
 and Dais[4] proved, in particular, 
 the equality of orbifold Hodge numbers and Hodge numbers of smooth crepant
 resolutions for Gorenstein global-quotient orbifolds. But there were 
 no results for non global-quotient orbifolds.\\
   
\noindent   In this paper, we show that the orbifold Hodge numbers of a generic 
 Calabi-Yau hypersurface in a complex 4-dimensional simplicial Fano       
 toric variety coincide with the Hodge numbers of its crepant resolution.
 Besides being the first non global-quotient example, this is also an
 important
 example in mirror symmetry. An immediate corollary of this is the mirror
 pairing of orbifold Hodge numbers for mirror families of Calabi-Yau
 3-folds in simplicial Fano toric varieties. It should be remarked here that
 the orbifold structure on the Calabi-Yau hypersurface that we consider, is the 
 one that naturally arises from its algebraic structure.\\  

\noindent  Recently the author has 
  observed that one can combine the characterisation of twisted sectors
  in this paper with the apparatus of string-theoretic Hodge numbers developed 
  by Batyrev and Dais [4], and Batyrev [2] to greatly generalise the results
  of this paper in the toric world. There is even a hope that the problem 
  may be completely solved, as far as Hodge numbers go, by using motivic 
  integration [7] and some stack theory.  But those results will be 
  discussed elsewhere.\\

\noindent Now we briefly describe how this article is organised. 
In sections 2 and 3 we 
review some facts from orbifold cohomology and toric geometry respectively.
In section 4 we find characterisations for the twisted sectors of simplicial
toric varieties and nondegenerate Calabi-Yau hypersurfaces of simplicial Fano
toric varieties. In section 5 we compute formulas for some orbifold Hodge 
numbers of these hypersurfaces, state some corollaries and then give an 
example. \\

\noindent\textbf{Acknowledgements}: This work was done under the 
supervision of my doctoral advisor,
Prof Yongbin Ruan. I am indebted to him for his guidance. My special 
thanks are also due to Professors Igor Dolgachev,  Lev Borisov and  David
Cox who were kind enough to listen to me and answer my questions. I thank 
Prof Cox for correcting an earlier draft of this article. I thank the faculty
and friends at UW-Madison for their help and kindness.

\section{Orbifolds}

\subsection{Orbifold Structure} 

 Let $U$ be a connected topological space, $V$ a connected $n$-dimensional
 smooth manifold and $G$ a finite group acting smoothly on $V$.
 An $n$-dimensional \emph{uniformising system} of $U$ is a triple 
 $(V,G,\pi)$,where $\pi:V \to U$ is a continuous map inducing a homeomorphism
 between the quotient space $V/G$ and $U$. Two uniformising systems $(V_i,
 G_i,\pi_i), i=1,2,$ are \emph{isomorphic} if there is a diffeomorphism 
 $\phi:V_1\to V_2$ and an isomorphism $\lambda:G_1 \to G_2$ such that $\phi$
 is $\lambda$-equivariant, and $\pi_2\circ \phi =\pi_1$.\\
 
 \textbf{Note:} If $(\phi,\lambda)$ is an automorphism of $(V,G,\pi)$, then
 there is a $g \in G$ such that $\phi(x)=g.x$ and $\lambda(a)=g.a.g^{-1}$,
 for any $x\in V$ and $a\in G$.\\
 
 Let $i:U'\to U$ be a connected open subset of $U$. An uniformising system
 $(V',G',\pi')$ of $U'$ is said to be induced from $(V,G,\pi)$ if there is a
 monomorphism $\lambda:G'\to G$ and a $\lambda$-equivariant open embedding
 $\phi:V'\to V$ such that $i\circ \pi' = \pi \circ \phi$. The pair $(\phi,
 \lambda):(V',G',\pi')\to(V,G,\pi)$ is called an \emph{injection}. Two 
 injections $(\phi_i,\lambda_i):(V'_i,G'_i,\pi'_i)\to(V,G,\pi), i=1,2,$ are
 \emph{isomorphic} if there is an isomorphism $(\psi,\tau)$ from $(V'_1,
 G'_1,\pi'_1)$ to $(V'_2,G'_2,\pi'_2),$ and an automorphism $(\bar{\psi},
 \bar{\tau})$ of $(V,G,\pi)$ such that $(\bar{\psi},\bar{\tau})\circ(\phi_1,
 \lambda_1) = (\phi_2,\lambda_2)\circ(\psi,\tau).$\\
 
 \textbf{Fact:} Let $(V,G,\pi)$ be an uniformising system of $U$. For any connected 
 open subset $U'$ of $U$, $(V,G,\pi)$ induces a unique isomorphism class of
 uniformising systems of $U'$.\\
 
 Let $U$ be a connected and locally connected topological space. For any point
 $p\in U$ we define the \emph{germ} of uniformising systems at $p$ in the
 following sense. Let  $(V_1,G_1,\pi_1)$ and $(V_2,G_2,\pi_2)$ be uniformizing
 systems of neighbourhoods $U_1$ and $U_2$ of $p$. We say that $(V_1,G_1,\pi_1)$  and $(V_2,G_2,\pi_2)$ are \emph{equivalent} at $p$ if they induce isomorphic
 uniformising systems for a neighbourhood $U_3$ of $p$.\\
 
\textbf{Definition}: Let $X$ be a Hausdorff, second countable topological
 space. An \emph{$n$-dimensional orbifold structure} on $X$ is given by 
 the following
 data : for any point $p\in X$, there is a neighbourhood $U_p$ and an 
 $n$-dimensional uniformising system $(V_p,G_p,\pi_p)$ of $U_p$ such that for 
 any point $q\in U_p$, $(V_p,G_p,\pi_p)$ and $(V_q,G_q,\pi_q)$ are equivalent
 at $q$.\\

\textbf{Definition}: The \emph{germ of an orbifold structure} on $X$ is an
 equivalence class under the following equivalence relation :
 two orbifold structures $\{(V_p,G_p,\pi_p):p\in X\}$ and 
 $\{(V'_p,G'_p,\pi'_p):p\in X\}$ are equivalent if for any $p\in X$,
 $(V_p,G_p,\pi_p)$ and $(V'_p,G'_p,\pi'_p)$ are equivalent at $p$. With a 
 given germ of orbifold structures on it, $X$ is called an \emph{orbifold}.\\

 We call each $U_p$ a uniformised neighbourhood of $p$, and $(V_p,G_p,\pi_p)$
 a chart at $p$. In fact we choose $U_p$ to be small enough that $G_p$ has the
 minimum possible order; that is, every element of $G_p$ fixes the preimage of
 $p$ in $V_p$. In what follows, this choice is assumed.\\ 
 
 An open subset $U$ of  $X$ is called an uniformised open set 
 if it is uniformised by $(V,G,\pi)$ such that for each $p\in U$,  $(V,G,\pi)$
 defines the same germ as $(V_p,G_p,\pi_p)$ at $p$. $X$ is called a \emph{global-quotient} orbifold if $X$ itself is an uniformised open set.
 A point $p$ is called 
 \emph{smooth} or \emph{regular} if $G_p$ is trivial; otherwise, it is called
 \emph{singular}. The set of smooth points is denoted by $X_{reg}$, and the set
 of singular points is denoted by $X_*$.\\

 An orbifold $X$ is called \textbf{reduced} if $G_p$ acts 
 effectively on $V_p$.
 We can canonically associate a reduced orbifold to an orbifold by redefining
 the group actions. Furthermore if a group element acts nontrivially, we require
 that the fixed-point set is of at least (real) codimension two, so that the 
 complement is locally connected. The toric varieties and hypersurfaces that
 we will deal with in this article are reduced orbifolds.\\ 
 
\textbf{Remark:}  It is important to note here that even a reduced non-smooth
 orbifold can have a smooth underlying variety because of examples with complex 
 reflections. For instance, consider $\cc/\zz_2$ where $\zz_2$ acts by 
 reflection about the origin. Then $\cc/\zz_2$ is complex analytically 
 isomorphic to $\cc$ but it has a singularity at the origin in the orbifold
 sense. \emph{Gorenstein} orbifolds do not present this problem as they do not
 admit such complex reflections.

\subsection{Orbifold Cohomology}

 First we will describe the so-called twisted sectors.
 Consider the set of pairs:\\

 $\tilde{X_*}=\{(p,(g)_{G_p})|p\in{X_*},1\neq g\in{G_p}\}$,\\

\noindent where $(g)_{G_p}$ denotes the conjugacy class of $g$ in $G_p$. Then 
 Kawasaki showed (see [5]) that $\tilde{X_*}$ is naturally an orbifold with an 
 orbifold structure given by \\

 $\{\pi_{p,g}:(V_p^g,C(g)/K_g)\to V_p^g/C(g) | p\in{X_*},1\neq{g}\in{G_p}\}$,\\

\noindent where $V_p^g$ is the fixed-point set of $g$ in $V_p, C(g)$ is the 
centralizer of $g$ in $G_p$, and $K_g$ is the kernel of the action of $C(g)$ 
on $V_p^g$. 
There is a map $\pi:\tilde{X_*}\to{X}$ defined by $(p,(g)_{G_P})\mapsto p$
such that $\pi(\tilde{X_*})=X_*$. In a sense $\tilde{X_*}$ is a resolution 
of singularities of $X$.\\

 Let $\tilde{X}=\tilde{X_*} \bigsqcup X$. We will describe the connected
 components of $\tilde{X}$. Recall that each point $p$ has a local chart  
 $(V_p,G_p,\pi_p)$ which gives a local uniformized neighbourhood 
 $U_p=\pi_p(V_p)$. If $q\in U_p$, up to conjugation, there is an injective
 homomorphism $G_q\to G_p$. For $g\in G_q$, the conjugacy class $(g)_{G_p}$
 is well-defined. We define an equivalence relation 
 $(g)_{G_q}\cong(g)_{G_p}$. Let $T$ denote the set of equivalence classes.
 By an abuse of notation, we use $(g)$ to denote the equivalence class to 
 which $(g)_{G_q}$ belongs. $\tilde{X}$ is decomposed as a disjoint union of 
 connected components\\

\noindent $\tilde{X}=\bigsqcup_{(g)\in T} X_{(g)},$\\
 where\\
 $X_{(g)}=\{(p,(g')_{G_p})|g'\in G_p,(g')_{G_p}\in (g)\}$.\\

 \textbf{Definition}: $X_{(g)}$ for $g\neq 1$ is called a \textbf{twisted sector}.
 We call $X_{(1)}=X$ the \emph{nontwisted} sector.

 Let $TX$ denote the \emph{orbifold tangent bundle} of $X$. See Chen and Ruan[5] for a definition. 
 An almost complex structure $J$ on $X$ is a smooth section of the orbifold 
 bundle $End(TX)$ such that $J^2=-Id.$ Assume that $X$ is an almost complex
 orbifold with an almost complex structure $J$. Then for any $p\in X_*$, $J$
 gives rise to an effective representation $\rho_p:G_p\to GL(n,\cc)$. For any
 $g\in G_p$ we write $\rho_p(g)$, upto conjugation, as a diagonal matrix \\
 $diag(e^{2\pi i\frac{m_{1,g}}{m_g}},....,e^{2\pi i\frac{m_{n,g}}{m_g}})$,
 where $m_g$ is the order of $g$ in $G_p$, and $0\leq m_{i,g}<m_g.$  Define 
 a function $\iota:\tilde X \to \qq$ by \\

 $ \iota(p,(g)_{G_p}) = \sum_{i=1}^n \frac{m_{i,g}}{m_g}.$\\

 This function $\iota:\tilde X \to \qq$ is locally constant. Denote it by
 $\iota_{(g)}$. It has the following properties:\\
 (1) $\iota_{(g)}$ is integral iff $\rho_p(g)\in SL(n,\cc)$.\\
 (2) $\iota_{(g)} + \iota_{(g^{-1})} = rank(\rho_p(g)-Id)=n-dim_{\cc}X_{(g)}$.

\noindent \textbf{Definition}: $\iota_{(g)}$ is called the \emph{degree shifting number} of $X_{(g)}$. \\

\noindent \textbf{Definition}: An almost complex orbifold is called 
 \emph{Gorenstein} if property(1) holds for all $(g)$.\\

\noindent \textbf{Remark:} An almost complex, complex or kahler structure on $X$ induces  a  corresponding similar structure on each $X_{(g)}$.\\

\noindent \textbf{Definition}: Let $\ff$ be any field containing $\qq$ as a subfield. 
 We define the \emph{orbifold chomology groups} of $X$ with coefficients in
 $\ff$ by \\

  $ H_{orb}^{d}(X;\ff) = \oplus_{(g)\in{T}} H^{d-2\iota_{(g)}}(X_{(g)};\ff)$.\\
 
\noindent \textbf{Definition}: Let $X$ be a closed complex orbifold. We define, for 
 $0 \le p,q \le dim_{\cc}X,$ \emph{ orbifold Dolbeault cohomology groups}\\

 $H_{orb}^{p,q}(X;\cc) = \oplus_{(g)}H^{p-\iota_{(g)},q-\iota_{(g)}}(X_{(g)};\cc)$.\\

\noindent \textbf{Remark:} When $X$ is a closed Kahler orbifold (so is each 
$X_{(g)}$),
these 
 Dolbeault groups are related to the singular cohomology groups of $X$ and
 $X_{(g)}$
 as in the manifold case, and the Hodge decomposition theorem holds for these 
 cohomology groups.\\

\noindent \textbf{Definition}: We define \emph{orbifold Hodge numbers} by 
 $h_{orb}^{p,q}(X)= dim H_{orb}^{p,q}(X;\cc)$.

\section{Facts from Toric Geometry}

\subsection{Fan, Orbits and Divisors}

 A complex $n$-dimensional toric variety is  constructed from  an
 $n$-dimensional lattice $N$ and a fan $\Xi$ in $N$. $\Xi$
 is a collection of strongly convex rational polyhedral cones 
 $\sigma$ in the real vector space $N_{\rr}=N\otimes_{\zz}\rr$, satisfying the 
 conditons: every face of a cone in $\Xi$ is also a cone in $\Xi$, and the 
 intersection of two cones in $\Xi$ is a face of each. We shall always assume
 that $\Xi$ is finite.\\

 Let $M=Hom(N,\zz)$ denote the dual lattice, with dual pairing denoted by 
 $\lan,\ran$.
 If $\sigma$ is a cone in $N$, the dual cone $\csig$ is the set of vectors in
 $M_\rr$ that are nonnegative on $\sigma$. This determines a commutative
 semigroup $\ssig = \csig \cap M =\{u\in M:\lan u,v\ran \geq 0$ for all $v\in
 \sigma\}.$
 This semigroup is finitely generated. The corresponding group-algebra
 $\cc[\ssig]$ which is the $\cc$-algebra with generators $\chi^m$ for each
 $m\in \ssig$ and relations $\chi^m\chi^{m'}=\chi^{m+m'}$, is a finitely 
 generated commutative $\cc$-algebra. Such an algebra corresponds to an
 $n$-dimensional affine variety $\Us:=$ spec$\cc[\ssig]$.\\
 
 If $\tau$ is a face of $\sigma$, then $\ssig$ is contained in $\stau$. So
 $\cc[\ssig]$  is a subalgebra of $\cc[\stau]$. This gives a map $\Ut \to \Us$.
 In fact, $\Ut$ is a principal open subset of $\Us$: if we choose $u\in \ssig$
 so that $\tau=\sigma \cap u^{\perp}$, then $\Ut = \{x\in \Us: u(x)\neq 0\}$.
 With these identifications, these affine varieties for different cones fit 
 together to  form an algebraic variety $X_\Xi$.  We will write $X$ for $X_\Xi$ when there is no confusion. $X$ is a normal
 Cohen-Macaulay variety.\\
 
 The properties of $\Xi$ strongly affect the geometry of $X$. For example:

\noindent (a)$X$ is complete(i.e., compact) iff the \emph{support}
 $|\Xi|:=\cup_{\sigma \in \Xi}\sigma = N_\rr$.\\
 (b)$X$ is nonsingular iff for every cone in $\Xi$, its 1-dimensional
 generators (defined below) are part of a $\zz$-basis of $N$.
 Such a fan is called \emph{smooth}.\\
 (c)$X$ is an orbifold iff the generators for every cone in $\Xi$ are linearly
  independent over $\rr$. We say $X$ and $\Xi$ are \emph{simplicial}.\\

 The affine variety corresponding to the trivial cone $\{0\}$ is the torus 
 $T_N=N\otimes \cc =$ Spec$\cc[M]$. It is an affine 
 (Zariski) open subset of $X_\Xi$.
 The action of $T_N$  on itself extends to an algebraic action of $T_N$ on
 $X_\Xi$. This action has exactly one orbit corresponding to each cone
 $\tau \in \Xi : O_{\tau} = Hom(\tau^{\perp} \cap M,\cc^*)$.
 We write $\tau<\sigma$ if $\tau$ is a face of $\sigma$ and we write 
 $\bar{O_{\sigma}}$ for the closure of $O_{\sigma}$ in $X$. Then \\

\noindent (a) $\Us = \bigsqcup_{\tau<\sigma}O_{\tau}$; [Note $\tau=\sigma$ is 
included]\\
 (b) $\bar{O_{\tau}} = \bigsqcup_{\gamma>\tau}O_{\gamma}$;\\
 (c) $X_\Xi = \bigsqcup_{\tau \in \Xi}O_{\tau}$;\\

 For each $d$, $\Xi(d)$ denotes the set of $d$ dimensional cones of $\Xi$.
 We reserve the letter $\eta$ to stand for elements of $\Xi(1)$. For each
 $\eta$, let $v_\eta$ denote the unique generator of the semigroup 
 $\eta \cap N$. Using these generators, a cone $\sigma \in \Xi$ can be 
 written $\sigma = \{ \sum_{\eta \subset \sigma} a_\eta v_\eta : a_\eta \ge 0\}$. The $v_\eta \in \sigma$ are the \emph{generators} of $\sigma$.
 A standard abuse of notation is to identify $\eta$ with
 $v_\eta$. If $r=|\Xi(1)|$ is the number of 1-dimensional cones, we sometimes
 write the $v_\eta$'s as $v_1,...,v_r$.\\

 Each $\bar{O_\eta}$ is an 
 irreducible $T_N$-invariant Weil divisor denoted $D_\eta$. Furthermore,
 $m \in M$ gives a character $\chi^m:T_N \to \cc^*$, and regarding $\chi^m$
 as a rational function on $X$, we have
 $div(\chi^m)=\sum_{\eta}\lan m,v_\eta \ran D_\eta $.
 We will always assume that the 1-dimensional cones span $N_\rr$. Then the
 Chow group $A_{n-1}(X)$ of Weil divisors modulo linear equivalence can be
 computed directly from the fan . By using the above , we get an exact 
 sequence \\

 $ 0\longrightarrow M\longrightarrow \zz^{\Xi(1)}\longrightarrow A_{n-1}(X)\longrightarrow 0$\\

\noindent where $m \in M$ maps to $(..,\lan m,v_\eta \ran,..)\in \zz^{\Xi(1)}$ 
and
$(..,a_\eta,..) \in \zz^{\Xi(1)}$ maps to the divisor class of 
$\sum_{\eta}a_{\eta}D_\eta$. Thus $A_{n-1}(X)$ has rank $r-n$.
 A Weil divisor $D=\sum_{\eta}a_{\eta}D_\eta$ is Cartier iff for each
 $\sigma \in \Xi$, there is $m_\sigma \in M$ such that
 $\lan m_\sigma,v_\eta \ran =-a_\eta$ whenever $\eta \subset \sigma$. 
 Polytopes arise naturally when dealing with toric varieties. 
 If $X$ is complete and $D=\sum_{\eta}a_{\eta}D_\eta$ is Cartier, then\\

\noindent $\Delta_D=\{m\in M_\rr:\lan m,v\ran \ge \phi_D(v) \forall v\in N_{\rr} \}\\
          =\{m\in M_\rr:\lan m,v\ran \ge -a_\eta \forall \eta \}$
	    is a polytope.\\
	    
 There is a $T_N$-equivariant map 
 $H^{0}(X,\mathcal{O}(D)) \simeq \bigoplus_{m\in \Delta_D \cap M} \cc \chi^m.$\\
 To see why this is true, think of sections of $\mathcal{O}(D)$ as rational 
 functions $f$ on $X$ such that div$(f) + D \ge 0$. Div$(\chi^m) +D\ge 0$ is
 equivalent to $\lan m,v_\eta \ran \ge -a_\eta \forall \eta$ and the isomorphism follows.\\

\noindent A cartier divisor $D$ is ample\\
 $\iff \lan m_\sigma,v_\eta \ran >-a_\eta$ whenever
 $\eta$ is not in $\sigma$ and $\sigma$ is $n$-dimensional. \\
 $\iff m_\sigma \neq m_\tau$ for $\sigma \neq \tau$ in $\Xi(n)$ and $\Delta_D$
 is an $n$-dimensional polytope with vertices $\{m_\sigma: \sigma \in \Xi(n)\}$.
 
\noindent A polytope is called integral if its vertices are integral.
$\Delta_D$ is integral if $D$ is ample.\\
 
\subsection{Projective Toric variety from polytope}

 Given any $n$-dimensional integral polytope $\Delta$ one can canonically
 associate a $\cc$-algebra $S_\Delta$ to it. See [[6],section 3.2.2]
 for details.
 Then the projective variety $\pp_\Delta:=$ Proj$(S_\Delta)$ turns out to be
 a toric variety with a fan called the \emph{normal fan} of $\Delta$. We will
 describe the normal fan in the special case of Fano toric varieties later.\\
 
 The $T_N$ orbit closures of $\pp_\Delta$ are in one-to-one
 correspondence with the nonempty faces $F$ of $\Delta$. There is a canonical
 surjection of polytope rings $\pi_F:S_\Delta \to S_F$ which induces a natural
 inclusion of toric varieties $\pp_F \hookrightarrow \pp_\Delta$.               
 $\pp_\Delta$ comes with a specific choice of ample divisor $D_\Delta$ such
 that $\Delta_{D_\Delta}= \Delta$. \\
  
  Choose a 
  basis for $M$. This corresponds to picking coordinates $t_1,...,t_n$ for the 
  torus $T_N$. Then, if $m \in M$ is written $m=(a_1,..,a_n)$, we have 
  $\chi^m = \prod_{i=1}^n t_{i}^{a_i}$, so we can write $t^m$ instead of 
  $\chi^m$.     
   For any $k\ge 0$, we have the space of Laurent polynomials 
 $ L(k\Delta)=\{f: f=\sum_{m\in k\Delta \cap M} \lambda_m t^m ,
 \lambda_m \in \cc \}.$
 Each $f\in L(k\Delta)$ gives the affine hypersurface $Z_f \subset T_N$ defined
 by $f=0$. 
 $L(k\Delta)\simeq H^{0}(\pp_\Delta,\mathcal{O}_{\pp_\Delta}(kD_\Delta)).$ 
 Under this  isomporphism, $f$ corresponds to an effective divisor
 $\bar{Z_f} \subset \pp_\Delta$. $\bar{Z_f}$ is a compactification of $Z_f$
 and is exactly the hypersurface of $\pp_\Delta$ corresponding to $f$.  
 We will use the following notation:\\
\noindent (a) $l(k\Delta) = |k\Delta \cap M | =$dim$(L(k\Delta))$ \\
 (b) $l^{*}(k\Delta)=| \{m \in k\Delta \cap M :m$ is not in any facet of
     $k\Delta \cap M \}|$\\
 
\subsection{Homogeneous coordinate ring}

 We introduced coordinates $t_1,..,t_n$
 on the torus of a toric variety. But it is useful to have global coordinates,
 similar to homogeneous coordinates on projective space.
  If $X$ is given by a fan $\Xi$ in $N_\rr$, we introduce a variable $x_\eta$
  for each $\eta \in \Xi(1)$ and consider the polynomial ring 
  $ S = \cc[x_\eta : \eta \in \Xi(1)].$ 
    A monomial in $S$ is written $x^D=\prod_{\eta}x_{\eta}^{a_\eta}$, where
  $D=\sum_{\eta}a_{\eta}D_\eta$ is an effective torus-invariant divisor on
  $X$. We say that $x^D$ has degree  $ deg(x^D)=[D] \in A_{n-1}(X)$.\\
  
  Thus, $S$ is graded by the Chow group $A_{n-1}(X)$. Given a divisor class 
  $\alpha \in A_{n-1}(X)$, $S_\alpha$ denotes the graded piece of $S$ of degree
  $\alpha$. We often write the variables as $x_1,...,x_r,$ where $x_i$
  corresponds to the cone generator $v_i$ and $r=|\Xi(1)|$. Then
  $S=\cc[x_1,...,x_r]$. The ring $S$, together with the 
  grading defined above is called 
  the \emph{homogeneous coordinate ring } of $X$. 
  We refer the reader to [[6],chapter 3.2] for a discussion on reconstructing
 the toric variety starting from the homogeneous coordinate ring.\\
 
 If $\tau$ is any cone of $\Xi$ then the orbit closure $\bar{O_\tau}$ is 
 given by the ideal $(x_i: i$ such that $v_i$ is a generator of $\tau)$ of
 $S$. Also
 the graded pieces of $S$ have nice cohomological interpretation. We noted that
   $L(\Delta)\simeq H^{0}(\pp,\mathcal{O}_{\pp}(D_\Delta)).$ Now the map 
   sending the Laurent monomial $t^m$ to $\prod_{\eta}x_{\eta}^{\lan m,v_\eta 
   \ran + a_\eta}$ induces an isomorphism 
   $H^{0}(X,\mathcal{O}_X (D)) \simeq S_\alpha$, where
   $\alpha = [D]\in A_{n-1}(X)$.\\

\subsection{Fano toric Varieties}
 
First we note that for an arbitrary toric variety $X$, the anticanonical divisor
 $-K_X = \sum_\eta D_\eta$. The following may be taken as definitions.
 
 \noindent \textbf{Definition}: A toric variety $X$ is
 \emph{Gorenstein} iff $-K_X$
 is Cartier. A complete toric variety $X$ is \emph{Fano} iff $-K_X$ is Cartier
 and ample.\\

 The anticanonical divisor of a Fano toric variety $X$ detemines an
 integral polytope
 $\Delta$ with special properties :
 
\noindent (a) All facets $\Gamma$ of $\Delta$ are supported by an affine hyperplane 
 of the form $\{m\in M_\rr:\lan m,v_\Gamma \ran = -1\}$ for some $v_\Gamma \in N$.\\
 (b)Int$(\Delta) \cap M = {0}$.

\noindent \textbf{Definition}: Such a polytope is called \emph{reflexive}.\\

  The \emph{polar polytope} $\Del$ of the reflexive polytope $\Delta$ is 
  obtained by $\Del = \{v \in N_\rr:\lan m,v \ran \ge -1$ for all
  $m\in \Delta\}\subset N_\rr$.
  The fan $\Xi$ of $X$ can be retrieved by coning over the proper faces of
  $\Del$. This fan is called the \emph{normal fan} of $\Delta$ and 
  $X= \pp_\Delta$.
  $\Del$ is also reflexive and $(\Del)^\circ = \Delta$. The Fano troic variety
  constructed from the normal fan of $\Del$ is denoted by $\pp_{\Del}$.\\

 We shall use $F$ and $F^\circ$ to denote a face of $\Delta$ and $\Del$
 respectively. There exists an inclusion reversing duality between the faces
 of $\Delta$ and $\Del$. For instance, the face of $\Delta$ which is dual to 
 the face $F^\circ$ of $\Del$ is defined as $\widehat{F^\circ} := \{m \in \Delta : \lan n,m \ran =-1 \forall n \in F^\circ \}$. If $\tau$ is the cone in the
 normal fan of $\Delta$ associated to the face $F^\circ$, then the orbit
 closure $\bar{O_\tau} = \pp_{\widehat{F^\circ}} $.\\
 
 Generic anticanonical hypersurfaces $V$ in 
 $\pp_{\Delta}$ and $V^{\circ}$ in $\pp_{\Del}$ constitute mirror 
 families of \textbf{Calabi-Yau} varieties. These varieties are orbifolds if 
 the corresponding ambient toric variety is simplicial. Observe that 
 $V$ can be identified with an element of $L(\Delta)$.\\
 
 A `maximal projective
 subdivision' of the fan of the toric variety gives a \emph{maximal projective
 crepant partial} (MPCP) resolution [see [1] or [6]]
 of the hypersurface by taking proper
 transforms. Let $\widehat{V}$ and $\widehat{V^{\circ}}$ denote
 the MPCP resolutions of $V$ and 
 $V^{\circ}$ respectively. These are again Calabi-Yau. These are smooth
 if $n=4$.\\
 
 \textbf{Remark:} A gorenstein toric variety or a Calabi-Yau hypersurface  
 of a Fano toric variety has gorenstein canonical singularities. In particular,
 the degree shifting numbers are integers.  Moreover,
 \emph{the singularity locus has at least complex codimension two}.  
 
\section{Twisted Sectors}

 We claim that the twisted sectors of a toric variety or a Calabi-Yau
 hypersurface
 can be identified with suborbifolds. Note that in general a twisted
 sector could be a multiple cover of the corresponding singular locus even
 if the group actions are all abelian.
 The situation is simpler in toric case because there exist canonical
 local coordinates  and uniformizing systems such that the group actions
 are in fact diagonal.
 They can be easily chased around to prove the claim. \\  

\subsection{Twisted Sectors in Simplicial Toric Variety}

 Let $\Xi$ be any simplicial fan. Then the orbifold structure of the
 toric variety $X_\Xi$ can be descibed as follows. 
 Let $\sigma$ be any n dimensional cone of $\Xi$.
 Let $v_1,...,v_n$  be the primitive 1  dim generators of
 $\sigma$. These are linearly independent in $N_\rr$.
 Let $\Ns$ be the sublattice of $N$ generated by $v_1,.., v_n$.
 Let $\Gs:= N/\Ns$ be the quotient group. $\Gs$ is finite and abelian.
 
 Let $\sigma'$ be the cone $\sigma$ regarded in $\Ns$. Let $\csigp$ be the 
 dual cone of $\sigma'$ in $\Ms$, the dual lattice of $\Ns$.
 $\Usp=$spec$(\cc[\csigp \cap \Ms])$. Note that $\sigma'$ is a smooth cone in
 $\Ns$. So $\Usp \cong \cc^n$.
 
 There is a canonical dual pairing $ \Ms /M \times N/\Ns \to \qq/\zz \to \cc^*$,
 the first map by the pairing $\lan ,\ran$ and the second by
 $q \mapsto$ exp$(2\pi iq)$. Now $\Gs$ acts on $\cc[\Ms]$ ,
 the group ring of $\Ms$, by : 
 $ v(\chi^u) =$ exp$(2\pi i\lan u,v\ran)\chi^u$, for $v \in N$ and $u\in \Ms$.
 Note that \\
 $(\cc [\Ms])^{\Gs} = \cc [M]$ \hfill .....(4.1.1) \\
 
  Thus $\Gs$ acts on $\Usp$. Let $\pi_\sigma$ be the quotient map.
 Then $\Us=\Usp/\Gs$. So 
 $\Us$ is uniformised by $(\Usp,\Gs, \Pi_\sigma)$. For any $\tau<\sigma$,
 the orbifold structure on $\Ut$ is same as the one induced from the
 uniformising system on $\Us$. Then by the description of the toric gluing 
 it is clear that $\{(\Usp,\Gs, \Pi_\sigma) :\sigma \in \Xi(n)\}$ defines a
 reduced orbifold structure on $X$. We give a more explicit verification of
 this fact below.\\
 
 \noindent Let $\B$ be the nonsingular matrix with generators $v_1,..,v_n$ of 
 $\sigma$ as rows.
 Then $\check{\sigma'}$ is generated in $\Ms$ by the the column vectors 
 $v^1,..,v^n$ of the matrix $\Bi$. 
 So $\chi^{v^1},..,\chi^{v^n}$ are the coordinates of $\Usp$.  
 For any $\kappa=(k_1,..,k_n)\in N$, the corresponding coset $[\kappa]\in \Gs $ 
 acts on $\Usp$ in these coordinates as a 
 diagonal matrix: diag$(\exp^{2\pi ic_1},...,\exp^{2\pi ic_n})$ where 
 $c_i = \lan \kappa,v^i \ran$.
 Such a matrix is uniquely represented by an $n$-tuple $a=(a_1,..,a_n)$
 where $a_i \in [0,1)$ and $c_i=a_i+b_i, b_i \in \zz$. In matrix notation,
 $\kappa \B^{-1}=a + b \iff \kappa = a\B+b\B.$
 We denote the inetgral
 vector $a\B$ in $N$ by $\kappa_a$ and the diagonal matrix corresponding
 to $a$ by $g_a$. ${\kappa_a} \leftrightarrow {g_a}$ gives a one to one
 correspondence between the elements of $\Gs$ and the integral vectors in 
 $N$ that are linear combinations of the generators of $\sigma$ with cofficients in $[0,1)$.\\

 Now let us examine the orbifold chart induced by $(\Usp,\Gs,\pi_\sigma)$
 at any point $x\in \Us$. By the orbit decomposition, there is unique 
 $\tau \in \sigma$ such that $x\in O_\tau$. Wlog assume $\tau$ is generated
 by $v_1,..,v_j ,j\le n$. Then any preimage of $x$ with respect to $\pi_\sigma$
 has coordinates $\chi^{v^i}=0$ iff $i \le j.$ Let $z=(0,..,0,z_{j+1},..,z_n)$
 be one such preimage. Let $\Gt:=\{g_a \in \Gs : a_i=0$ if $j+1\le i\le n\}$.   
 We can find a small neighbourhood 
 $W\subset (\cc^*)^{n-j}$ of $(z_{j+1},..,z_n)
 $ such that the inclusions 
 $\cc^j\times W\hookrightarrow \Usp$ and $\Gt \hookrightarrow \Gs$ induces an
 injection of uniformising systems $(\cc^j\times W,\Gt,\pi) \hookrightarrow 
 (\Usp,\Gs,\pi_\sigma)$ on some small open neighbourhood $U_x$ of $x$. So
 we have $G_x = \Gt$ and an orbifold chart
 $(\cc^j \times W,\Gt,\pi)$.
 Note that $\Gt$ can be constructed from the set 
 $\{\kappa_a=\sum_{i=1}^j a_iV_i : \kappa_a \in N, a_i \in [0,1)\}$ which 
 is completely determined by $\tau$ and hence is independent of $\sigma$.\\

  Now we determine the twisted sectors. Take any $x \in X$. $x$ belongs to a 
  unique $O_\tau$. Wlog the generators of $\tau$ are $v_1,..,v_j$. Consider 
  any $n$-dimensional $\sigma > \tau$. Wlog $v_1,..,v_n$ generates $\sigma$.
  Pick $g_a$ in $G_x$ such that $a_i \neq 0, \forall i\le j $ i.e,
  $\kappa_a$ lies in the interior of $\tau$. We want to find the twisted
  sector
  $X_{(g_a)}$. Consider $g_a$ as an element of $\Gs$. It is clear that
  $g_a$ fixes 
  $z \in \Usp$ iff $z_1= .. =z_{j+s}=0,$ for some $s\ge 0$ i.e.,
  $\pi_\sigma(z) \in O_\tau$ or $\pi_\sigma(z) \in  O_\delta$ for some
  $\delta>\tau$.
  So $X_{(g_a)} \cap \Us = \bar{O_\tau}\cap \Us$. Since a twisted sector
  is connected,  $X_{(g_a)}=\bar{O_\tau}$.
  If $g_{a} \in G_x$ is such that (wlog) only $ a_1,..,a_k \neq 0, k<j$, then
  $g_{a} \in G_\delta$ where $\delta$ is the cone generated by $v_1,..,v_k$ 
  and by the above argument $X_{(g_a)}=\bar{O_\delta}$. Thus we have proved the  following theorem.
 
\begin{theorem} A twisted sector of any simplicial toric variety $X_\Xi$
 is isomorphic to a subvariety $\bar{O_\tau}$ of $X_\Xi$ for some cone 
 $\tau \in \Xi$. Moreover, there is a one-to-one correspondence between the
 set of twisted sectors of the type $\bar{O_\tau}$ and the set of integral 
 vectors in the interior of $\tau$ which are linear combinations of the 
 1-dimensional generators of $\tau$ with coeffients in $(0,1)$.
\end{theorem} 

 Note that the degree shifting number $\iota_{(g_a)} = \sum a_i$.
 Now if $X_\Xi$ is Fano, i.e., $\Xi$  is obtained by coning over the 
 faces of a reflexive polytope $\Del$, then the twisted sectors with
 $\iota = 1$ are in one to one correspondence with the integral 
 interior points of faces of $\Del$. \\
 
\subsection{Twisted sectors of hypersurface of a Fano variety}

 We identify the twisted sectors of a generic nondegenerate anticanonical
 (Calabi-Yau) hypersurface $V$ of a simplicial Fano toric variety 
 $X=\pp_\Delta$. \textbf{Nondegenerate} means that
 $V\cap O_\tau$ is either empty or a smooth 
 subvariety of codimension one in $O_\tau$, for each torus
 orbit $O_\tau$ in $X$. Then $V$ turns out to be a suborbifold of $X$.
 Also nondegeneracy is a generic condition. For a different treatment on
 these, see [3].
 We show that $V=\bar{Z_f}$, for a generic $f \in L(\Delta)$, is nondegenerate
 and a suborbifold of $X$.\\ 
 
 Consider any $n$-dimensional cone $\sigma$
 with generators $v_1,...,v_n$. For 
 notational simplicity set $\chi^{v^i} = z_i$. Then $z_1,..,z_n$ are the 
 coordinates of $\Usp$. Let $Y$ be the preimage of $V\cap \Us$ in $\Usp$ with
 respect to $\pi_{\sigma}$. Then $Y$ is defined by the equation 
 $\sum_{m \in \Delta \cap M}\lambda_m \prod_{i=1}^n {z_i}^{\lan m,v_i\ran +1} = 0$. This is because, $t^m =\prod_{i=1}^n {z_i}^{q_i} \iff m=\sum q_i v^i =\Bi q
 \iff q = m\B  \iff q_i= \lan m,v_i \ran$. The one is added  
 to ensure that $V$ is anticanonical. Note that by definition of $\Delta$,
 $\lan m,v_i \ran + 1 \ge 0$. If $\lambda_{m_\sigma} \ne 0$ then $Y$ does 
 not pass through the origin.
 It can be checked from this description using Bertini's theorem
 that for generic values of the 
 coefficients $\lambda_m, Y$ is a smooth submanifold of $\Usp$ 
 that intersects
 the coordinate planes $z_{i_1}=..=z_{i_j}=0$ transversely. \\

 $Y$ is $\Gs-$stable by (4.1.1). When $Y$ is smooth, all singularities of
 $V\cap \Us$ are quotient singularities induced by action of $\Gs$ on $Y$.
 Since there are only
 finitely many $n$-dimensional cones, $V$ is nondegenerate and a suborbifold of
 $X$.
 $(Y,\Gs,\pi_\sigma)$ is an uniformizing system for $V \cap \Us$. 
 Let $\tau$ be the face of $\sigma$ obtained by coning over the face $F^{\circ}$ of
 $\Del$. Wlog let $v_1,..,v_j : j\ge 2$ be the
 generators of $\tau$. (We remarked earlier in section 3.4 that the Gorenstein
 condition rules out codimension one singularities.)  
 We want to find a chart for any point $x\in V\cap O\tau$.  By our
 earlier remark that $Y$ misses the origin of $\Usp$,
 $V\cap \bar{O_\sigma}$ is empty. So we need only consider proper faces of 
 $\sigma$.
 By a result of Fulton[[8],section 5.3], $V \cap \bar{O_\tau}$ consists 
 of $l^*(\widehat{F^{\circ}})+1$ points if $F^\circ$ is a codimension $2$ face
 of $\Del$ i.e., $j=n-1 \iff$
 dim$O_\tau =1$. Since the the only other points in $\bar{O_\tau}$ in this case  are $O_\sigma$ for $n$-dimensional cones $\sigma >\tau$, all the intersection
 points actually lie in $O_\tau$.
 If codimension $F^\circ$ is bigger than $2$, then
 again by Bertini $V \cap \bar{O_\tau}$ is irreducible.\\

 Following the same notation as before, $x$ has a small neighborhood 
 $U_x \cap Y$  such that $((\cc^j \times W)\cap Y,\Gt,\pi)$ is a chart for
 $V$ at $x$. $\cc^j \times W$ is , as before, a suitable neighbourhood of 
 some preimage $z$ of $x$ in $\Usp$. The tangent space $TY_z$ is a 
 $\Gt$-stable subspace of $T\cc^{n}_z$. Any $g_a \in \Gt$ acts trivially 
 on $TW_z =$span$\{\partial/{\partial z_i}, i=j+1,..,n\}.$ By transversality,
 we can choose basis $\{\xi_1,..,\xi_n \}$ of $T\cc^{n}_z$ such that $\xi_i \in
 TY_z$ $\forall i\le n-1$ and $\xi_n \in TW_z$. $g_a$ acts trivially on $\xi_n$.
 This implies that the degree shifting number of $g_{a}|_{TY_z}$ is still
 $\sum_{i=1}^j a_i$.\\

 From the description of the charts, it is clear that twisted sectors of $V$
 are isomorphic to $V\cap \bar{O_\tau}$ where $2 \le$dim$(\tau) \le n-1.$
 Recall that $\bar{O_\tau}= \pp_{\widehat{F^\circ}}$ where ${\widehat{F^\circ}}$
 is the face of $\Delta$ dual to $F^\circ$.
 In particular we have the following theorem.

\begin{theorem} Let $V$ be a generic nondegenerate anticanonical hypersurface
 of an $n$-dimensional simplicial Fano toric variety $\pp_{\Delta}$. Then the
 twisted sectors of $V$ are isomorphic to $V\cap \pp_{\widehat{F^\circ}}$ for 
 some face $F^\circ$ of $\Del$ such that $1 \le$dim$F^\circ \le n-2$. There is
 exactly one
 twisted sector of this type having $\iota =1$, corresponding to each integral 
 interior point of
 $F^\circ$ if dim$F^\circ < n-2$. If dim$F^\circ = n-2$, then there are exactly
 $l^*(\widehat{F^\circ})+1$ twisted sectors of this type having  $\iota =1$,
 corresponding to each integral interior point of $F^\circ$.
\end{theorem}

\section{Orbifold Hodge Numbers}

\subsection{$h^{1,1}_{orb}(V)$} 

 Let $V_{(g)}$ denote a twisted sector of the hypersurface $V$
 and $\iota_{(g)}$ its degree shifting number.

\noindent $h^{1,1}_{orb}(V) = h^{1,1}(V) + \sum_{\iota_{(g)} =1} h^{0,0}(V_{(g)})$.
 Since $ h^{0,0}(V_{(g)}) =1 $ for each twisted sector, by theorem (4.2)  we
 obtain\\
 
 $\sum_{\iota_{(g)} =1} h^{0,0}(V_{(g)})$\\
  $= \sum_{1\le dim(F^\circ) \le n-2} l^*(F^\circ) + \sum_{dim(F^\circ)=n-2}
  l^*(F^\circ)l^*(\widehat{F^\circ})$\\
  $=l(\Del)-r-1-\sum_{dim(F^\circ)=n-1}l^*(F^\circ) + \sum_{dim(F^\circ)=n-2}
  l^*(F^\circ)l^*(\widehat{F^\circ})$\\
 
\noindent To compute $h^{1,1}(V)$ we invoke the following Lefschetz hyperplane
 theorem [[3],Proposition 10.8] .
 
\begin{lemma} Let $V$ be a nondegenerate
 ample hypersurface of an $n$-dimensional complete simplicial toric variety
 $X$. Then the natural map,induced by inclusion,
 $j^{*}: H^i (X) \to H^i (V)$ is an isomorphism for 
 $i<n-1$ and an injection for $i=n-1$.
\end{lemma}

\noindent In our case $V$ is anticanonical, and since the anticanonical 
 divisor of a Fano variety is ample, $V$ is ample. Also it is well 
 known[[8],section 5.1]
 that for any simplicial toric variety $X$, $H^{2}(X,\rr)=H^{1,1}(X,\rr)=
 A_{n-1}(X)\otimes \rr$.
 So for $n \ge 4$ , 
 $h^{1,1}(V)=h^{1,1}(\pp_\Delta)=$rank$A_{n-1}(\pp_\Delta)=r-n$. 
 Thus we have the following theorem.

\begin{theorem} For any generic nondegenerate anticanonical 
hypersurface $V$ of an
 $n$-dimensional simplicial Fano toric variety $\pp_\Delta$, $n \ge 4$,\\ 
 
 $h^{1,1}_{orb}(V) =l(\Del)-n-1-\sum_{dim(F^\circ)=n-1}l^*(F^\circ) +
 \sum_{dim(F^\circ)=n-2}l^*(F^\circ)l^*(\widehat{F^\circ})$. 
\end{theorem}
 
 \subsection{$h^{n-2,1}(V)$}
 
 Next we compute $h^{n-2,1}(V)$ for $n \ge 4$.
 For this we have to use the homogeneous 
 coordinate ring $S$ of $X=\pp_\Delta$. Let $v_1,..,v_r$ be the one dimensional
 cones of the normal fan of $\Delta$. Let $x_1,..,x_r$ be the corresponding
 homogeneous coordinates. 
 Let $\beta_{0}=[-K_X]=[\sum_{i=1}^r D_i] \in A_{n-1}(X)$.
 Then $S_{\beta_0} \simeq L(\Delta)$. And the divisor $f\in L(\Delta)$
 corresponding to $V$ can be written in the homogeneous coordinates as 
 $ f = \sum_{m\in M\cap L(\Delta)} \lambda_m \prod_{i=1}^r {x_i}^{\lan m,v_i\ran +1}.$
 For notational simplicity, we will denote 
 $\prod_{i=1}^r {x_i}^{\lan m,v_i\ran +1}$ by $\mathbf{x}^m$ for any $m\in M$.
 Define the Jacobian ideal of $f$ to be 
 $J(f)=\lan \partial{f}/\partial{x_1},...,\partial{f}/\partial{x_r} \ran.$\\
 
 First we quote the following theorem [[3], theorem 10.13] .

 \begin{lemma} Let $X$ be an d-dimensional complete
 simplicial toric variety and $V\subset X$ be a quasi-smooth(i.e.,suborbifold)
 ample hypersurface defined by $f\in S_\beta$.  Then for $k\ne (d/2)+1$, there
 exists a canonical isomorphism

 $ (S/J(f))_{k\beta -\beta_0} \simeq PH^{d-k,k-1}(V)$.\\ 
 \end{lemma}

\noindent\textbf{Remark:} The primitive cohomolgy
$PH^{d-1}(V):=H^{d-1}(V)/($im$H^{d-1}(X).$
 This coincides with the usual cohomology if $d$ is even.\\

 For our application of this lemma: $d=n, k=2$ and $\beta = \beta_0$. So we need   to compute the rank of $(S/J(f))_{\beta_0}$. Also, $H^{n-2,1}(X) =0$ if 
 $n \ge 4$ , so that $H^{n-2,1}(V)=PH^{n-2,1}(V)$.\\
 
\begin{lemma}
 $x_i \partial{f}/\partial{x_i} \in (J(f))_{\beta_0}, i=1,..,r,$ and the 
 space of complex linear relations among these has dimension $r-(n+1)$.
\end{lemma}

\noindent\textbf{Proof}

 $x_i \partial{f}/\partial{x_i} = \sum_m \lambda_m (\mvi +1)\mathbf{x}^m$.
 
 $\sum_i c_i x_i\partial{f}/\partial{x_i} =  \sum_m \lambda_m
 (\lan m,\sum_i c_i v_i \ran + \sum_i c_i)\mathbf{x}^m$\\
 
\noindent For a generic $f$ we can assume that $\lambda_m \ne 0$ for each $m\in \Delta\cap M$. Hence $\sum_i c_i x_i\partial{f}/\partial{x_i} \equiv 0 \iff 
 \lan m,\sum_i c_i v_i \ran + \sum_i c_i = 0$ $\forall m\in \Delta \cap M $

\noindent  In particular, taking $m=0$ we get $\sum_i c_i =0$. Therefore 
 $\lan m,\sum_i c_i v_i \ran = 0$ $\forall m\in \Delta \cap M $ and since 
 $\Delta$ is $n$-dimensional we have $\sum_i c_i v_i = 0$.

\noindent  Thus $\sum_i c_i x_i\partial{f}/\partial{x_i}\equiv 0 \iff \sum_i c_i v_i = 0$
 and $\sum_i c_i =0$.\\

\noindent  Now let $\tilde{v_i}=(v_i,1)\in \rr^{n+1}\subset \cc^{n+1}=\rr^{n+1}\otimes \cc
$.
 
\noindent  Since the $v_i$'s are vertices of $\Del$, the $\tilde{v_i}$'s  are generators
 of the $(n+1)$-dimensional cone $\{ q\in \rr^{n+1}: 
 qt\in \Del$ for some $ t\in \rr_{> 0} \}$.
 So the $\tilde{v_i}$'s span $\rr^{n+1}$ over $\rr$, and hence they span 
 $\rr^{n+1}\otimes \cc$ over $\cc$.

\noindent  Note that $(\sum c_i v_i, \sum c_i) = \sum c_i \tilde{v_i}$.
 
\noindent  Hence the lemma follows.\\

 Wlog assume that the $ x_k\partial{f}/\partial{x_k}, k=1,..,n+1,$ are linearly
 inedendent. In other words $\tilde{v_k}, k=1,...,n+1$ are linearly independent.
 We want to find monomials $\prod_{j\ne i}{x_j}^{p_j}$ that have same degree
 in $S$ as $x_i$. So we want $m^*\in M$ such that $\lan m^*,v_j\ran = p_j \ge 0$
 if $j\ne i$ and $\lan m^*,v_i\ran =-1$. Such $m^*$ is given by interior lattice
 points of $F_i$, the $(n-1)$-dimensional face of $\Delta$ that is dual to  the $0$-dimensional face $\{v_i\}$ of $\Del$.\\
 
 Then for each $m^*\in $(Int$F_i \cap M$), 
 $\prod_{j\ne i}{x_j}^{\msvj} \partial{f}/\partial{x_i}$ belongs to
 $(J(f))_{\beta_0}$. Together with the $x_i \partial{f}/\partial{x_i}$, 
 these generate $(J(f))_{\beta_0}$ as we vary over all $i$.\\
 
\noindent $\prod_{j\ne i}{x_j}^{\msvj} \partial{f}/\partial{x_i}$\\

\noindent  $= \prod_{j=1}^r {x_j}^{\msvj} x_i\partial{f}/\partial{x_i}$\\
 
\noindent  $=\sum_{m\in \Delta \cap M} \lambda_m (\mvi +1)\mathbf{x}^{m+m^*}$\\
 
\noindent  $=\sum_{m'\in \Delta \cap M \cap(\Delta +m^*)} \lambda_{m'-m^*} 
 (\lan m'-m^*,v_i\ran +1)\mathbf{x}^{m'}$\\
 
\noindent  $=\sum_{m\in \Delta \cap M} \lambda_{m-m^*}(\lan m-m^*,v_i\ran +1)
 I(m-m^* \in \Delta)\mathbf{x}^m $\\

 Let Int$F_i \cap M =\{m_{i,i_s}: 1\le s\le t_i ; t_i \ge 0\}$.
 
 Then $(J(f))_{\beta_0} =$span$\{x_k\partial{f}/\partial{x_k}, 
  \prod_{j=1}^r {x_j}^{\lan m_{i_s},v_j\ran} x_i\partial{f}/\partial{x_i}
  : 1\le k \le n+1, 1\le i_s \le t_i , i=1,..,r \}.$ 
 We want to find  the dimension of this complex vector space.
 So we study the space of linear relations :\\

\noindent $\sum_k c_k x_k\partial{f}/\partial{x_k} + \sum_{i,i_s} d_{i,i_s}
 \prod_{j=1}^r {x_j}^{\lan m_{i_s},v_j\ran} x_i\partial{f}/\partial{x_i}
 \equiv 0 $

\noindent  $\iff$ 
 
\noindent  $\sum_m \{ \sum_k c_k \lambda_m (\mvk +1) + \sum_{i,i_s} d_{i,i_s}
 \lambda_{m-m_{i,i_s}} I(m-m_{i,i_s}\in \Delta) (\lan m-m_{i,i_s},v_i\ran +1) \}
 \mathbf{x}^m \equiv 0$ 

\noindent  $\iff$

\noindent  $ \sum_k c_k \lambda_m (\mvk +1)    + \sum d_{i,i_s}
 \lambda_{m-m_{i,i_s}} I(m-m_{i,i_s}\in \Delta) \mvi \equiv 0 $
 for each $m\in \Delta \cap M$,\hfill [note: $\lan m_{i,i_s},v_i\ran =-1$] \\

 This is a system of $l(\Delta)$ number of linear equations in 
 $\gamma =n+1+\sum_{i=1}^r l^* (F_i)$  variables namely $ c_k , d_{i,i_s}$.
 Note that $l(\Delta)\ge \gamma$. We shall find 
 a nonsingular subsystem of rank $\gamma$.\\
 
 To do 
 so pick $n$  linearly independent vertices $m_1,...,m_n$ of $\Delta$
 and let $m_{n+1} = 0$, the origin. Then from the above system we pick
 the equations corresponding to $m = m_1,...,m_{n+1}$ and $m= m_{i,i_s}:
 i=1,..,r; 0\le i_s \le t_i$. Denote this $\gamma \times \gamma$ system by (**).
 It can be written as :
 $$
   \left[
   \begin{array}{ll}
   \mathbf{P} & \mathbf{A} \\
    \mathbf{B} & \mathbf{Q}
   \end{array}
   \right]
   \left(
   \begin{array}{ll}
   c \\
   d
   \end{array}
   \right)
   =
   \left(
   \begin{array}{ll}
   0 \\
   0
   \end{array}
   \right)
$$ 
\noindent where
 $$
   \mathbf{P}=
   \left[
   \begin{array}{lll}
   \lambda_{m_1}(\lan m_1,v_1\ran +1) & ... & \lambda_{m_1}(\lan m_1,v_{n+1}\ran +1)\\
   ... & ... & ...\\
   \lambda_{m_n}(\lan m_n,v_1\ran +1) & ... & \lambda_{m_n}(\lan m_n,v_{n+1}\ran +1)\\
   \lambda_0 & ... & \lambda_0
   \end{array}
   \right]
$$

$$
 \mathbf{Q}=
 \left[
 \begin{array}{lll}
 \lambda_{m_{1,1}-m{1,1}}I(.)(\lan m_{1,1}-m_{1,1},v_1\ran +1) 
 & ... & \lambda_{m_{1,1}-m_{r,t_r}}I(.)(\lan m_{1,1}-m_{r,t_r},v_r\ran +1) \\
  ... & ... & ... \\
  \lambda_{m_{r,t_r}-m{1,1}}I(.)(\lan m_{r,t_r}-m_{1,1},v_1 \ran +1)
 & ... & \lambda_{m_{r,t_r}-m_{r,t_r}}I(.)(\lan m_{r,t_r}-m_{r,t_r},v_r\ran +1) 
 \end{array}
 \right]
 $$\\

 Observe that all the diagonal entries of $\mathbf{Q}$ are $\lambda_0$,
 and none of its off-diagonal entries has $\lambda_0$. Also any entry of
 $\mathbf{A}$ is of the form $\lambda_{m_k-m_{i,i_s}}I(.)$ and hence does
 not involve $\lambda_0$. Similarly an entry of $\mathbf{B}$ is of the form
 $\lambda_{m_{i,i_s}}(\lan m_{i,i_s},v_k \ran +1)$ and so does not have 
 $\lambda_0$.\\
 
 Consider the  determinant of the coefficient matrix 
$
  \left[
  \begin{array}{ll}
   \mathbf{P} & \mathbf{A} \\
    \mathbf{B} & \mathbf{Q}
   \end{array}
   \right]
$
as a polynomial in the $\lambda$'s. Then the term of this determinant
having the highest power of $\lambda_0$ is 
$(\lambda_0)^{\sum l^*(F_i)}$det$\mathbf{P}$. We will show below that 
det$\mathbf{P}$ = nonzero constant times $\lambda_{m_1}...\lambda_{m_n} \lambda_0$.
Thus the determinant of the coeffcient matrix of the system  (**)
is a nontrivial polynomial in the $\lambda$'s and is therefore nonzero for 
generic choice of the $\lambda$'s.  Hence $(J(f))_{\beta_0}$ has rank $\gamma$
as a complex vector space, for a generic $f \in L(\Delta).$
Since $S_{\beta_0} \simeq L(\Delta)$, so $(S/J(f))_{\beta_0}$ has rank
$l(\Delta)- \gamma$, for a generic $f$.\\
  
\begin{lemma} The $(n+1)\times (n+1)$ matrix $\mathbf{P}
 =((P_{i,j}=\lambda_{m_i}(\lan m_i,v_j\ran +1) ))$ is 
 nonsingular for generic chioce of $\lambda$'s.
\end{lemma}

\noindent \textbf{Proof}
 
 Let $\mathbf{E}$ be the $(n+1)\times(n+1)$ matrix
 $((E_{i,j}=(\lan m_i,v_j\ran +1) ))$.  Then 
 det$\mathbf{P}= \lambda_{m_1}...\lambda_{m_{n+1}}$det$\mathbf{E}$.

 We claim that $\mathbf{E}$ is nonsingular. Otherwise there exists a nontrivial
 vector $(c_1,...,c_{n+1})$ such that $\sum_{k=1}^{n+1} c_k (\lan m_i,v_k\ran +1 ) =0$ for all $i=1,..,n+1$. In particular, for $i=n+1$ we get $\sum_{k=1}^{n+1} c_k =0$. This implies  $\sum_{k=1}^{n+1} c_k \lan m_i,v_k\ran =0$ for all
 $i=1,..,n$. Since $m_1,..,m_n$ are linearly independent, this 
 would imply that $\sum_{k=1}^{n+1} c_k \lan m,v_k\ran =0$ for all
 $m\in \Delta$. Therefore $\sum_{k=1}^{n+1}c_k v_k =0$. this combined with
 $\sum c_k =0$ implies that $\sum_{k=1}^{n+1} c_k \tilde{v_k} =0$ which
 contradicts the linear independence of $\tilde{v_1},..,\tilde{v_{n+1}}$.
 Thus the lemma holds.\\
 
 So we have the following theorem.
\begin{theorem} For any generic nondegenerate anticanonical hypersurface $V$ of an $n$-dimensional simplicial Fano toric variety $\pp_\Delta$, and
$n \ge 4 $,\\

 $h^{n-2,1}(V)= l(\Delta)-n-1-\sum_{dim(F)=n-1} l^*(F)$.
\end{theorem}

\noindent Note : $F$ denotes a facet of the polytope $\Delta$ in the
statement of the above theorem.

\subsection{Cohomology of the twisted sectors of $V$}

\noindent We now want to compute $h^{n-3,0}(V_{(g)})$ for any twisted sector
 $V_{(g)} \cong V\cap \pp_{\widehat{F^\circ}}$. This is obviously zero
 if dim$(F^\circ) > 1$, since dim $(\pp_{\widehat{F^\circ}})= n-1-$dim$(F^\circ)$. So we will only consider the case dim$(F^\circ) = 1$. Let $\tau$ be the 
 $2$-dimensional cone obtained by coning over $F^\circ$. As noted earlier
 $\bar{O_\tau}= \pp_{\widehat{F^\circ}}$.
 The restriction of $V$ to $\bar{O_\tau}$ gives a quasismooth  ample 
  hyprsurface of  $\bar{O_\tau}$, which we shall identify with
  $V_{(g)}$. So we are again in a situation where we can invoke lemma 2.\\

  For this we need to understand the homogeneous coordinate ring $S'$ of 
 $\bar{O_\tau}$. According to Fulton[[8],section 3.1], a fan for $\bar{O_\tau}$ can be
 constucted from the fan $\Xi$ of $X$ as follows.\\
 
 Let $N_\tau$ be the sublattice of $N$ generated by the primitive one
 dimensional generators of $\tau$. Let $N(\tau) =N/N_\tau$.
 The dual lattice of $N(\tau)$ is given by $M(\tau)={\tau}^\perp \cap M$.
 The \emph{star} of a cone $\tau$ can be defined abstractly as the set of cones
 $\sigma$ in $\Xi$ that contain $\tau$ as a face. Such cones $\sigma$ are 
 determined by their images in $N(\tau)$ i.e. by 
$\bar{\sigma}= (\sigma + (N_\tau)_{\rr}) /(N_\tau)_{\rr} 
\subset N_{\rr}/(N_\tau)_{\rr} =N(\tau)_{\rr}.$ 
 These cones $\{\bar{\sigma}:\tau < \sigma \}$ form a fan Star$(\tau)$ in 
 $N(\tau)$. $\bar{O_\tau}$ is the toric variety corresponding to this fan.
 
 Wlog let $v_1,v_2$ be the generators of $\tau$. The
 corresponding Weil divisors in $X$ are $D_1$ and $D_2$, and 
 $\bar{O_\tau}$ is the intersection 
 $|G_\tau| D_1 D_2$ considered as an element of the
 Chow ring of $X$.
 
 Again wlog assume that
 $v_j, j=3,..,l$ are the one dimensional cones of $\Xi$
 such that $\{v_1,v_2,v_j\}$
 generate a three dimensional cone of $\Xi$. In other words,
 $\bar{v_j}, j=3,..,l$
 are the one dimensional cones of Star$(\tau)$. Let $\tilde{D_j}:=D_1 D_2 D_j$
 for $j=3,..,r$. Note that $\tilde{D_j} = 0$ if $j>l$. 
 The divisor of  $\bar{O_\tau}$ corresponding to $\bar{v_j}$ is $\tilde{D_j}$
 for $j=3,..,l$. 
 So the homogeneous coordinate ring $S'$
 is generated by variables $y_j$ corresponding to $\tilde{D_j}$ for $j=3,..,l$.
 Denote by $\alpha_0$ the anticanonical class $\sum_{j=3}^l \tilde{D_j}$ in $S'$.\\
 
 On the other hand the divisor $V$ restricts to $-K_X D_1 D_2 
 =(D_1+...+D_r)D_1D_2 = \sum_{j=3}^l \tilde{D_j} +(D_1+D_2)D_1D_2$.
 To see what $(D_1+D_2)D_1D_2$ is in terms of the $\tilde{D_j}$s , we can
 pick a point $m \in \widehat{F^\circ}\cap M$ and let $b_i:=\lan m,v_i\ran ,
 1\le i\le r$. Then $\sum_{i=1}^r b_iD_i$ is linearly equivalent to zero.
 Note that $b_1=b_2=-1$. Hence we have $D_1+D_2 =\sum_{i=3}^r b_iD_i$. So, 
 $(D_1+D_2)D_1D_2= \sum_{i=3}^l b_i \tilde{D_i}$.
 Let $\alpha$ be the class in $S'$ representing the effective ample divisor
 $-K_X D_1 D_2$. Let $f'$ be the associated homogeneous polynomial in the 
 $y_j$s. Now we can apply lemma 2 to the $(n-2)$-dimensional variety 
 $\bar{O_\tau}$ and the ample hypersurface $V_{(g)}$. Choose $k=1$ in the
 lemma to get \\
 $$(S'/J(f'))_{\alpha -\alpha_0} \simeq PH^{n-3,0}(V_{(g)})=H^{n-3,0}(V_{(g)})$$
 since $ H^{n-3,0}(\bar{O_\tau})= 0$.

 Now $\alpha -\alpha_0 =[\alpha -\sum_{j=3}^l \tilde{D_j}]$.
 A typical generator $\partial{f'}/\partial{y_i} \in J(f')$ has degree
 $[\alpha - \tilde{D_i}]$. There are no nonconstant
 regular functions on the projective variety $\bar{O_\tau}$.
 So any nontrivial effective divisor, and in particular $\tilde{D_i}$, 
 $\sum_{j=3,..l;j\ne i} \tilde{D_j}$
 and $\sum_{j=3,..l}\tilde{D_j}$ are not linearly equivalent to zero.
 This implies that $(J(f'))_{\alpha -\alpha_0} =0$. Hence we obtain that 
 $$(S')_{\alpha -\alpha_0} \simeq H^{n-3,0}(V_{(g)})$$

 Now  $\alpha -\alpha_0 =[\sum_{j=3}^l b_j\tilde{D_j}]$. We want to 
 identify the effective divisors in this class. So we want $m_* \in M(\tau)$
 such that $\sum_{j=3}^l (b_j+\lan m_*,\bar{v_j} \ran) \tilde{D_j}$ is 
 effective.  This is if and only if $(b_j+\lan m_*,\bar{v_j}\ran)\ge 0$
 for all $j=3,..,l$\\ $\iff  (\lan m+m_*,v_j\ran)\ge 0 $ for $j=3,..,l$ \\
 $\iff m+m_* \in $Int$\widehat{F^\circ}\cap M$. 
 
 To justify the last step note that $\lan m+m_*,v_i\ran =-1$ for $i=1,2$
 so that $m+m_* \in \widehat{F^\circ}\cap M$. If $m+m_*$ is not in the 
 interoir of $ \widehat{F^\circ}$ then it belongs to some face 
 $ \widehat{E^\circ}$ where $F^\circ \subset E^\circ$ as faces of $\Del$.
 The cone $\delta$ associated to $E^\circ$ contains $\tau$ as a face and
 hence contains some $v_k$ for $3\le k\le l$. Then $\lan m+m_*,v_k \ran =-1$
 which is a contradiction.

 Since $m$ is fixed, the required effective divisors are in one-to-one
 correspondence with the interior lattice points of $ \widehat{F^\circ}$.
 Hence $h^{n-3,0}(V_{(g)})= l^*(\widehat{F^\circ})$. Since there are 
 $l^*(F^\circ)$ twisted sectors isomorphic to $V\cap \pp_{\widehat{F^\circ}}$
 we have the following.\\
 
\noindent $h^{n-2,1}_{orb}(V)$\\
 $=h^{n-2,1}(V) + \sum_{\iota(g)=1} h^{n-3,0}(V_{(g)})$\\
 $=l(\Delta)-n-1- \sum_{dim(F)=n-1}l^*(F) + \sum_{dim(F^\circ)=1}l^*(F^\circ)
 l^*(\widehat{F^\circ})$\\
 $=l(\Delta)-n-1- \sum_{dim(F)=n-1}l^*(F) + \sum_{dim(F)=n-2}l^*(F)
 l^*(\widehat{F}). $\\
 
 For the last step we used the one-to-one correspondence between faces of 
 $\Delta$ and $\Del$.

\subsection{Main Results}

 Combining the above formula with theorem 4 we have the following theorem.

\begin{theorem} For any generic nondegenerate anticanonical hypersurface $V$ of an n-dimensional simplicial Fano toric variety $\pp_\Delta$, $n\ge 4$, \\

 $ h^{n-2,1}_{orb}(V)= l(\Delta)-n-1- \sum_{dim(F)=n-1}l^*(F)+\sum_{dim(F)=n-2} l^*(F)l^*(\widehat{F})$.
 
\end{theorem}

\begin{corollary} If $\widehat{V}$ is an MPCP desingularisation of any generic   nondegenerate anticanonical hypersurface $V$ of an
 n-dimensional simplicial Fano toric variety $\pp_{\Delta}$, $n\ge 4,$ then
 $h^{p,1}_{orb}(V)=h^{p,1}(\widehat{V})$ for $p=1$ and $p=n-2$
\end{corollary}

\noindent \textbf{Proof} The formulas for $h^{p,1}(\widehat{V})$ for $p=1,n-2$
computed
 in [1] by Batyrev match the orbifold Hodge numbers for $V$ obtained in
 theorem 3 and theorem 5.

\begin{corollary} In the case $n=4$, $h^{p,q}_{orb}(V)=h^{p,q}(\widehat{V})$
 for any $p$ and $q$.
\end{corollary}

\noindent \textbf{Proof} We need  only consider $p,q \le 3$.
 Also $h^{p,0}_{orb}\equiv h^{p,0}$ by definition since $\iota$ is nonnegative.
 So, by Serre duality, it is 
enough to consider just the cases $p=1,q=1$ and $p=2,q=1$. These are addressed
by corollary 1.

 We should remark here that in this case $\widehat{V}$ is actually smooth.

\begin{corollary} If $\pp_{\Del}$ is also simplicial, and $V^\circ$ is a 
 generic nondegenerate anticanonical hypersurface of $\pp_{\Del}$ , then
 $h^{1,1}_{orb}(V)=h^{n-2,1}_{orb}(V^\circ)$ and vice versa.
\end{corollary}

\noindent \textbf{Proof} Follows from interchanging the roles of
$\Delta$ and $\Del$ in the formulas.\\
 
\noindent \textbf{Remark} In particular, for the $n=4$ case, we have 
$h^{p,q}_{orb}(V)=h^{3-p,q}_{orb}(V^\circ)$. This is an example of `mirror
symmetry' of orbifold hodge numbers.
 
 \subsection{An Example}

 Consider the complex 4-dimensional weighted projective space 
 $X = \pp$(1,1,2,2,2). It is a simplicial Fano toric variety. It's fan
 $\Xi$ has the following 1-dimensional cones in $N \cong \zz^4$:
 $ v_1 = (-1,-2,-2,-2)$, $v_2 = (1,0,0,0)$, $v_3 = (0,1,0,0)$,
  $ v_4 = (0,0,1,0)$, $v_5 = (0,0,0,1). $
  $\Xi$ has five 4-dimensional cones, obtained by dropping one of the 
  $v_i$'s at a time and taking the cone generated by the remaining four.\\
    
  Let $D_i$ denote the torus-invariant divisor given by the orbit 
  closure $\bar{O_{v_i}}$. It is easy to check that in $A_3(X)$,
  $[D_2]=[D_1]$ and $[D_i]= 2[D_1]$ for $i \ge 3$. Construct the homogeneous 
  coordinate ring of $X$ by introducing variables $x_i$ corresponding to $v_i$.
  Then deg($x_1$)= deg($x_2$) = 1 (= $[D_1]$) and deg($x_i$)=2 for $i \ge 3$.\\

  This leads to the more familiar description of $\pp$(1,1,2,2,2) as 
  $(\cc^5 - \{0\})/\cc^*$. The action of any $\alpha \in \cc^*$ on 
  $\cc^5 - \{0\}$ is as follows :\\
 $ \alpha.[x_1:x_2:x_3:x_4:x_5] = [\alpha x_1:\alpha x_2: \alpha^2 x_3:\alpha^2 x_4:\alpha^2 x_5] .$\\

  In this description, $D_i$ corresponds to the hyperplane $\{x_i=0\}$ and
  the 4-dimensional cones of $\Xi$ correspond to the open sets $\{x_i \neq 0\}$.
  It is also easily seen that the singular locus of $X$ is precisely the 
  surface $\{x_1=x_2=0 \}$. In fact, this represents the only twisted sector
  of $X$. The $v_i$'s are the vertices of a reflexive polytope $\Del$.
  The faces of $\Del$ have only one interior lattice point : 
  (0,-1,-1,-1) = 1/2 $(v_1 + v_2)$. This lattice point corresponds to
  the twisted sector and the local isotropy group is $\zz_2$.\\ 
  
  The dual reflexive polytope $\Delta$ in $M_\rr$  has the following vertices:\\
  $ w_1=(-1,-1,-1,-1)$, $ w_2=(7,-1,-1,-1)$, $w_3 =(-1,3,-1,-1)$,\\
   $w_4=(-1,-1,3,-1)$, $w_5=(-1,-1,-1,3)$.\\

   $\Delta$ is the polytope corresponding to the anticanonical divisor 
   $-K_X = \sum_{i=1}^5 D_i$ of $X$, and $X = \pp_{\Delta}$. If $V$ is a
   generic nondegenerate Calabi-Yau (anticanonical) hypersurface of $X$, 
   then $V$ has just one twisted sector namely $V \cap \{x_1=x_2=0 \}$.
   One can directly compute the genus of this curve by using the 
   Riemann-Hurwitz formula. It turns out to be 3.\\

   The following is a list of the Hodge numbers of $V$:                       %
   $ h^{1,0} = h^{2,0} = 0$, $h^{3,0} = 1$, $h^{1,1} = 1$, $h^{2,1} = 83$,
  $ h^{1,1}_{orb} = 2$,  $h^{2,1}_{orb} = 86$.

   The dual Fano veriety $\pp_{\Del}$ is also simplicial. This is easily
   checked since its fan is is obtained by coning over the faces of $\Delta$.
   In fact, $\pp_{\Del}$ = $\pp_{\Delta}/ \zz_4^3$. This is also shown easily.
   First, observe that $w_1 = -w_2 - 2w_3 -2w_4 -2w_5 $. Secondly, if $\bar{M}$
   is the sublattice of $M$ generated by $w_2,w_3,w_4,w_5$, then 
   $M/\bar{M}$ = $\zz_4^3$.

   The mirror Calabi-Yau family is also easy to write down. Let $z_i$ denote 
   the homogeneous coordinate of $\pp_{\Del}$ correspoding to $w_i$. Then a 
   generic Calabi-Yau hypersurface $V^{\circ}$ of $\pp_{\Del}$ has equation:
   
   $\lambda_1 z_1^8 + \lambda_2 z_2^8 + \lambda_3 z_3^4 + \lambda_4 z_4^4 
   + \lambda_5 z_5^4 + \lambda_6 z_1^4 z_2^4 + \lambda_0 z_1z_2z_3z_4z_5 = 0.$\\

   Note that $h^{1,1}_{orb}(V^{\circ}) = 86$ and $h^{2,1}_{orb}(V^{\circ}) = 2.$ 
\section*{Bibliography}

\noindent  [1] V.Batyrev, \emph{Dual polyhedra and Mirror symmetry for
 Calabi-Yau hypersurfaces in toric varieties}, AG/9310003\\

\noindent  [2] V.Batyrev, \emph{Non-Archimedean Integrals and Stringy Euler
 numbers of log-terminal pairs}, JEMS 1, 1999.\\

\noindent  [3] V.Batyrev and D.Cox, \emph{On the Hodge structure of projective 
 hypersurfaces in toric varieties}, AG/9306011\\

\noindent  [4] V.Batyrev and D.Dais, \emph{Strong McKay correspondence,
 string-theoretic Hodge numbers and mirror symmetry}, AG/9410001 \\

\noindent  [5] W.Chen and Y.Ruan, \emph{A New cohomology theory for Orbifold},
 AG/0004129\\

\noindent  [6] D.Cox and S.Katz, \emph{Mirror Symmetry and Algebraic Geometry}, 
 Mathematical Surveys and Monographs, Vol 68, AMS.\\

\noindent  [7] J.Denef and F.Loeser, \emph{Germs of arcs on singular algebraic
 varieties and motivic integration}, Invent. Math. \textbf{135} (1999).\\
 
\noindent  [8] W.Fulton, \emph{Introduction to toric varieties},
 Princeton University Press, 1993.\\

\noindent  [9] M.Reid, \emph{La correspondance de McKay}, AG/9911165 \\

\noindent  [10] Y.Ruan, \emph{Stringy geometry and Topology of Orbifolds}, 
     AG/0011149  \\

\end{document}